\documentclass{amsart}
\numberwithin{equation}{section}
\usepackage{psfrag}
\usepackage{graphicx}
\usepackage[all]{xy}
\usepackage{amscd}
\usepackage{amssymb}
\let\cal\mathcal

\def\Lscr{{\cal L}}

\def\Oscr{{\cal O}}

\let\blb\mathbb
\def\CC{{\blb C}}

\def \PP{{\blb P}}

\def \NN{{\blb N}}

\def \HH{{\blb H}}

\def\Id{\operatorname{id}}
\def\pr{\mathop{\text{pr}}\nolimits}

\def\Der{\operatorname{Der}}
\def\Bimod{\operatorname{Bimod}}

\def\Mod{\operatorname{Mod}}
\def\mod{\operatorname{mod}}

\def\Lie{\operatorname{Lie}}

\def\Rep{\operatorname {Rep}}

\def\Hom{\operatorname {Hom}}

\def\tr{\operatorname {tr}}

\def\ker{\operatorname {ker}}

\def\r{\rightarrow}

\newtheorem{lemma}{Lemma}[section]
\newtheorem{proposition}[lemma]{Proposition}
\newtheorem{theorem}[lemma]{Theorem}

\newtheorem{convention}[lemma]{Convention}

\newtheorem{lemmas}{Lemma}[subsection]
\newtheorem{propositions}[lemmas]{Proposition}
\newtheorem{theorems}[lemmas]{Theorem}
\newtheorem{corollarys}[lemmas]{Corollary}

\theoremstyle{definition}

\newtheorem{example}[lemma]{Example}

\newtheorem{definitions}[lemmas]{Definition}

{

}

\theoremstyle{remark}

\newtheorem{remark}[lemma]{Remark}

\newdimen\uboxsep \uboxsep=1ex
\def\uboxn#1{\vtop to 0pt{\hrule height 0pt depth 0pt\vskip\uboxsep
\hbox to 0pt{\hss #1\hss}\vss}}

\def\uboxs#1{\vbox to 0pt{\vss\hbox to 0pt{\hss #1\hss}
\vskip\uboxsep\hrule height 0pt depth 0pt}}

\def\Gl{\operatorname{Gl}}

\def\DR{\operatorname{DR}}
\def\HH{\operatorname{HH}}
\def\HC{\operatorname{HC}}

\def\ldb{\mathopen{\{\!\!\{}} \def\rdb{\mathclose{\}\!\!\}}}

\def\Ad{\operatorname{Ad}}
\begin{document}
\author{Michel Van den Bergh}
\address{Departement WINI, Universiteit Hasselt, 3090 Diepenbeek, Belgium}
\email{michel.vandenbergh@uhasselt.be}
\thanks{The author is a director of Research at the FWO}
\title{Non-commutative quasi-Hamiltonian spaces}
\begin{abstract}
  In this paper we introduce non-commutative analogues for the
  quasi-Hamiltonian $G$-spaces introduced by Alekseev, Malkin and
  Meinrenken. We outline the connection with the non-commutative
  analogues of quasi-Poisson algebras which the author had introduced earlier.
\end{abstract}
\maketitle
\tableofcontents
\section{Introduction}
There has been recent interest in developing a non-commutative version
of differential geometry based on Kontsevich's philosophy
\cite{KoRo,Leb4} that for a property of a non-commutative $k$-algebra
$A$ to have geometric meaning it should induce standard geometric
properties on all representation spaces $\Rep(A,N)=\Hom(A,M_N(k))$.
Non-commutative symplectic geometry was developed in
\cite{LebBock,CBEG,Ginzburg1,Kosymp} and a non-commutative version of
(quasi-)Poisson geometry was introduced in \cite{VdB33}.

In this paper we introduce so-called \emph{quasi-bisymplectic algebras} (see
\S\ref{ref-6-32}). These are a multiplicative analogue of algebras equipped
with a bisymplectic form (see \cite{CBEG}). Our definition is such
that the representation spaces of quasi-bisymplectic algebras are
quasi-Hamiltonian $G$-spaces, as introduced in \cite{AMM}. We develop
non-commutative analogues for some aspects of the commutative theory
\cite{AKM,AMM}.  In particular we show that there is a one-one
correspondence between quasi-bisymplectic algebras and Hamiltonian
double quasi-Poisson algebras (introduced in \cite{VdB33}) which satisfy a
suitable non-degeneracy condition (see Theorem \ref{ref-7.1-35} below).

As a side result we show that double quasi-Poisson algebras give rise
to something we call a ``double Lie algebroid'' (see Theorem
\ref{ref-5.3-29} below). This is a non-commutative version of \cite[Thm.\
2.5]{BuCr}.

In the final section of the paper we show that the Hamiltonian double
quasi-Poisson algebras derived from quivers which were introduced in
\cite{VdB33} are non-degenerate. Hence these algebras are also in a
natural way quasi-bisymplectic.

The main result of this paper was announced in \cite[App A]{VdB33} where
we discussed the relation between \emph{ordinary} double Poisson brackets and
bisymplectic forms (i.e.\ the ``non-quasi''case). It should be said
however that our proof for the equivalence between  integrability
of double quasi-Poisson brackets and quasi-bisymplectic forms is
based on a brute force computation and is less satisfactory than the
corresponding proof in \cite[App A]{VdB33}.

This paper depends rather heavily on \cite{CBEG,VdB33}. For the
convenience of the reader we have included some preliminary sections
explaining the relevant concepts and results. To simplify the exposition
we have chosen to write out all our computations over a base ring
which is a field, although that is not sufficient for the application
to quivers. Therefore in the short section \S\ref{ref-8.1-57} we outline
the modifications necessary to handle more general situations. 

A change in presentation with respect to \cite{VdB33} is that
throughout the paper we have emphasized a certain functor
\[
(-)_N:\Bimod(A)\r \Mod(\Oscr(\Rep(A,N)))
\]
which connect an algebra with its $N$'th representation space. When
applied to a non-commutative object this functor yields the
corresponding classical object. For example if $L$ is a double Lie
algebroid over $A$ (see above) then $L_N$ is a classical Lie algebroid
on $\Rep(A.N)$.

The author wished to thank Victor Ginzburg for explaining some aspects
of~\cite{Ginzburg2}.
\section{Preliminaries}
\subsection{Representation spaces}
\label{ref-2.1-0}
We assume that $k$ is a field of characteristic zero although this
hypotheses is often too strong. 
Throughout $A$ is a
finitely generated $k$-algebra.
For $N\in \NN$ the associated representation space of $A$ is defined
as
\[
\Rep(A,N)=\Hom(A,M_N(k))
\]
The group $\Gl_{{N}}$ acts on $\Rep(A,{{N}})$ by conjugation on $M_{{N}}(k)$.

A natural point of view in non-commutative algebraic
geometry is that for
a property of a non-commutative ring $A$ to have geometric meaning
it should induce standard geometric properties on all $\Rep(A,{{N}})$. 

It is easy to see that $\Rep(A,{{N}})$ is an affine variety and its
coordinate ring has a very convenient description. It is easy to see
that the ring $A_N\overset{\text{def}}{=}\Oscr(\Rep(A,{{N}}))$ is generated by the symbols
$(a_{ij})_{ij=1,\ldots,N}$, subject to the relations
\[
(ab)_{ij}=a_{il}b_{lj}
\]
together with additivity in $a$ and $1_{ij}=\delta_{ij}$. Here and
below we sum over repeated indices.  While this description is very
convenient it is of course also very uneconomical. For example if
$A=k\langle (x_l)_{l=1,\ldots,n}\rangle$ then
$\Oscr(\Rep(A,N))=k[(x_{l,ij})_{l=1,\ldots,n,i,j=1,\ldots,N}]$.  It is
easy to verify that our assumption that $A$ is finitely generated
implies that $\Oscr(\Rep(A,N))$ is finitely generated.

If $a\in A$ then $a_{ij}$ defines a matrix values function on $\Rep(A,N)$ which
we sometimes denote by $X(a)$.
Concretely this is the function which associates to every $\phi:A\r M_N(k)\in
\Rep(A,N)$ the matrix $(\phi(a))_{ij}$. 

We define the trace $\tr(a)$ of $a\in A$ as $a_{ii}$. This defines a
$\Gl_n$-invariant function on $\Rep(A,{{N}})$. As $\tr([A,A])=0$
we see that elements of $A/[A,A]$ correspond to invariant functions
on $\Rep(A,N)$. In fact this will be true in all cases we
consider below: the non-commutative version of working with $\Gl_N$ invariants
objects is working modulo commutators. 
\subsection{Differential forms}
We now consider differential forms. The  bimodule $\Omega_A$ of
differentials is generated as an $A$-bimodule by the symbols $da$ subject
to the relations $d(ab)=a(db)+(da) b$ and linearity. 
As usual one puts $\Omega A=T_A \Omega_A$.
Defining $d(da)=0$, $d(a)=da$ makes $\Omega_A$ into a 
differential graded algebra. 
Every homogeneous element $\omega$ in $\Omega A$ has a representations $a_0
da_1\cdots da_n$. To such an element one associate a matrix valued
differential form:
\begin{equation}
\label{ref-2.1-1}
\omega_{ij}=a_{1,i_1}da_{2,i_1i_2}\cdots da_{n,i_{n-1}j}
\end{equation}
$(\omega_{ij})_{ij}$ is a matrix valued differential form on $\Rep(Q,\alpha)$. 
If we write it as $X(\omega)$ then \eqref{ref-2.1-1} may be rewritten as
\[
X(\omega)=X(a_1)dX(a_2)\cdots dX(a_n)
\]
It would be tempting to define the non-commutative de Rham complex of
$A$ as $\Omega A$. However, quite remarkably (see \cite[\S 2.5]{CBEG}), $\Omega A$ is acyclic. That is
\[
H^m(\Omega A)=
\begin{cases}
k &\text{$m=0$}\\
0&\text{otherwise}
\end{cases}
\]
Nevertheless $\Omega A$ can be used as the basis for a new 
construction of the cyclic homology of $A$ (see \cite{Ginzburg2}). 

A different non-commutative analogue of the de Rham complex is
the \emph{Karoubi-de Rham complex} which is defined by
\[
\DR(A)=\Omega A/[\Omega A,\Omega A]
\]
The de Rham cohomology of $A$ is defined as the cohomology of the
Karoubi complex. It is closely related to cyclic homology and to
equivariant de Rham cohomology of representation spaces (see \cite[Thm
4.2.3]{Ginzburg2}).  According to \cite[2.6.7]{Loday1} we have a short
exact sequence of reduced (co)homology
\[
0\r \bar{H}^n(\DR(A))\r \overline{\HC}_n(A) \r \overline{\HH}_{n+1}(A)\r 0
\]
Below we will mostly deal with \emph{smooth algebras}, i.e. algebras whose
category of bimodules has homological dimension one. In that case
$\overline{\HH}_{n}(A)=0$ for $n>1$ and hence $\bar{H}^n(\DR(A))=
\overline{\HC}_n(A)$ for $n\ge 1$.  The following lemma is instructive.

To any element $\omega$ we associate a
$\Gl_n$ invariant differential form  $\tr(\omega)=\omega_{ii}$. 
It this way we obtain a map
\[
\tr:DR(A)\r \Omega(\Rep(A,N))^{\Gl(n)}
\]
which descends to cohomology. 
\begin{example} \label{ref-2.1-2} Let $C=k[t,t^{-1}]$. 
Then
\[
H^i(\DR(C))=
\begin{cases}
k&\text{if $i=0$}\\
k(t^{-1}dt)^i&\text{if $i$ is odd}\\
0&\text{otherwise}
\end{cases}
\]
That the cohomology groups have the indicated dimension
follows from \cite[Cor 3.4.15]{Loday1}. It is easy to see that the
generators for the cohomology groups are as claimed. 

We have $\Rep(A,N)=\Gl_N$. The elements $\tr(t^{-1}dt)^{2i+1}$ are precisely
the generators of the de Rham cohomology of $\Gl_N$. 
\end{example}
\subsection{Vector fields}
\label{ref-2.3-3}
Now we discuss vector fields.  Again there are two possible points of view.

 If we insist
that a vector field on $A$  induces vector fields on all $\Rep(A,{{N}})$ 
then a vector field on $A$ should simply be a derivation $\Delta:A\r A$. 
The induced derivation $\delta$ on $\Oscr(\Rep(A,{{N}}))$ is then given by
\[
\delta(a_{ij})=\Delta(a)_{ij}
\]
A second point of view is that a vector field $\Delta$ on $A$ should induce
\emph{matrix valued vector fields} $(\Delta_{ij})_{i,j=1,\ldots,n}$ on all $\Rep(A,{{N}})$. Since now 
$\Delta_{ij}(a_{uv})$ depends on four indices $\Delta(a)$ should be an element
of $A\otimes A$. It was first of observed by Crawley-Boevey that the
second point of view is often more useful. Put  
\[
D_{A}\overset{\text{def}}{=}\Der(A,A\otimes A)=\Hom_{A\otimes A^\circ}
(\Omega_A,A\otimes A)
\] where 
as usual we put the outer bimodule structure on $A\otimes A$.  The 
corresponding matrix valued vector fields on  $\Rep(A,{{N}})$ are then given
by
\begin{equation}
\label{ref-2.2-4}
\Delta_{ij}(a_{uv})=\Delta(a)'_{uj}\Delta(a)''_{iv}
\end{equation}
where by convention we write an element $x$ of $A\otimes A$ as
$x'\otimes x''$ (i.e.\ we drop the summation sign). We will
call the peculiar arrangement of indices in \eqref{ref-2.2-4}
the \emph{standard index convention}. It will reappear often below.

 Starting with
$D_A$ we define \emph{the algebra of p{{oly-vector}} fields} $DA$ on
$A$ as the tensor algebra $T_A D_A$ of $D_A$ where we make $D_A$ into
an $A$-bimodule by using the inner bimodule structure on $A\otimes A$.
Any homogeneous element $\delta$ of $DA$ induces polyvector fields
$X(\delta)$ on all representation spaces using a formula similar to
\eqref{ref-2.1-1}. 

The counterpart to the differential on $\Omega A$ is the ``double
Schouten Nijenhuys'' bracket on $DA$ which was defined in \cite{VdB33}.
A double bracket on an ordinary algebra $A$
is a bilinear map 
\[
\ldb-,-\rdb:A\times A\r A\otimes A
\]
which is a derivation in its second argument (for the outer bimodule
structure on $A$) and which satisfies
\[
\ldb a,b\rdb =-\ldb b,a\rdb^\circ
\]
where $(u\otimes v)^\circ=v\otimes u$. If $\ldb-,-\rdb$ satisfies the following
analogue of the Jacobi identity
\[
0=\ldb a,b,c\rdb\overset{\text{def}}{=}\ldb a,\ldb b,c\rdb \rdb _L
+\tau_{(123)}\ldb b,\ldb c,a\rdb \rdb _L
+\tau_{(132)}\ldb c,\ldb a,b\rdb \rdb _L
\]
where for $\tau\in S_n$ we define
\begin{equation}
\label{ref-2.3-5}
\tau(a_1\otimes \cdots \otimes a_n)=a_{\tau^{-1}(1)}\otimes \cdots \otimes 
a_{\tau^{-1}(n)}
\end{equation}
(with sign in the graded case) then we call $\ldb-,-\rdb$ a double
Poisson bracket.  A double Poisson bracket induces a Lie bracket
$\{-,-\}$ on $A/[A,A]$ via the formula
\[
\{a,b\}=\ldb a,b\rdb'\ldb a,b\rdb''
\]
The motivation for introducing double Poisson brackets
is that they induce  ordinary Poisson brackets on $\Rep(A,N)$ via the standard index convention. The precise
formula is
\[
\{ a_{ij},b_{uv}\}=\ldb a,b\rdb'_{uj} \ldb a,b\rdb''_{iv}
\]
One of the main results of \cite{VdB33} is the following.
\begin{proposition}  The graded algebra
$DA$ has the structure of a \emph{double Gerstenhaber algebra}
  i.e.\ a (super) double Poisson algebra with a double Poisson bracket
$\ldb-,-\rdb$ of degree $-1$. 
 \end{proposition}
For the convenience of the reader we give the construction of the
double Schouten-Nijenhuys bracket on $D_A$. 

If $\delta,\Delta\in D_{A}$ then  it is easy to see that
\begin{align*}
  \ldb \delta,\Delta\rdb \,\tilde{}_l&=(\delta\otimes 1)\Delta-(1\otimes \Delta)\delta\\
  \ldb \delta,\Delta\rdb \,\tilde{}_r&=(1\otimes
  \delta)\Delta-(\Delta\otimes 1)\delta=-\ldb \Delta,\delta\rdb\,\tilde{}_l
\end{align*}
define derivations $A\r A^{\otimes 3}$ for the outer bimodule structure
on $A^{\otimes 3}$.
Since $\Omega_{A}$ is finitely generated we obtain
\[
\Der_B(A,A^{\otimes 3})\cong\Hom_{A^e}(\Omega_{A/B},A\otimes A)\otimes A
\]
We view $\ldb \delta,\Delta\rdb \,\tilde{}_l$ and $\ldb
\delta,\Delta\rdb \,\tilde{}_r$ as elements of $D_{A}\otimes_k A$
and $A\otimes_k D_{A}$ respectively.  To this end we define
\begin{align*}
\ldb \delta,\Delta\rdb _l&=\tau_{(23)}\circ \ldb \delta,\Delta\rdb \,\tilde{}_l\\
\ldb \delta,\Delta\rdb _r&=\tau_{(12)}\circ \ldb \delta,\Delta\rdb \,\tilde{}_r
\end{align*}
and we write
\[
\begin{split}
\label{ref-2.3-6}
\ldb \delta,\Delta\rdb _l&=\ldb \delta,\Delta\rdb '_l\otimes \ldb \delta,\Delta\rdb ''_l\\
\ldb \delta,\Delta\rdb _r&=\ldb \delta,\Delta\rdb '_r\otimes \ldb \delta,\Delta\rdb ''_r
\end{split}
\]
with $\ldb \delta,\Delta\rdb ''_l,\ldb \delta,\Delta\rdb '_r\in A$,
$\ldb \delta,\Delta\rdb '_l,\ldb \delta,\Delta\rdb ''_r$ in $D_{A}$.

An easy verification shows that
\[
\ldb \delta,\Delta \rdb_r=-\ldb \Delta,\delta \rdb_l^\circ
\]

The double Schouten-Nijenhuys bracket is defined on generators by
\[
\label{ref-2.3-7}
\begin{split}
\ldb a,b\rdb &=0\\
\ldb \delta,a\rdb &=\delta(a)\\
\ldb \delta,\Delta\rdb &=\ldb \delta,\Delta\rdb _l+\ldb \delta,\Delta\rdb _r
\end{split}
\]
for $a,b\in A$, $\delta,\Delta\in D_{A}$,
where we regard the righthand sides in the previous display as elements
of $DA\otimes DA$.

One may show again that the double Schouten-Nijenhuys bracket on $D_A$
induces the standard Schouten-Nijenhuys bracket on the algebra of
polyvector fields on $\Rep(A,N)$ using the standard index convention.

In \cite{VdB33} it was shown how to associate a double bracket $\ldb-,-\rdb_P$
to an element $P\in D^2A$. The formula is obtained by linear extension
from the following formula with $\delta,\Delta\in D_A$:
\begin{equation}
\label{ref-2.4-8}
\ldb a,b\rdb_{\delta\Delta}= \Delta(b)'\delta(a)''\otimes \delta(a)'\Delta(b)''
-\delta(b)'\Delta(a)''\otimes \Delta(a)' \delta(b)''
\end{equation}
If $A$ is smooth then the double Jacobi identity for $\ldb-,-\rdb$ is
equivalent to $\{P,P\}=0$.

The algebra $DA$ has a remarkable element $E$ which has no commutative
counter part. It is the double derivation which sends $a$ to $a\otimes
1-1\otimes a$. It is intimately connected to the $\Gl_N$ action
on $\Rep(A,N)$. More precisely the following result was proved in \cite{CBEG,VdB33}
\begin{propositions}
\label{ref-2.3.1-9}
  Let $f_{ij}\in M_{\alpha}=\Lie(\Gl_N)$ be the elementary matrix which is
  $1$ in the $(i,j)$-entry and zero everywhere else. 
Then $(E_p)_{ij}$ acts as $f_{ji}$ on $\Oscr(\Rep(A,N))$.
\end{propositions}
The element $E$ appears in several pleasing formulas. For example for any
poly-vector field $\Delta$ in $DA$ we have
\[
\ldb E,\Delta \rdb=\Delta\otimes 1-1\otimes \Delta
\]
A vector field of the form
\[
H_a=\ldb a,-\rdb
\]
is called \emph{Hamiltonian}. An element $\Phi\in A$ is a \emph{moment map} 
if it realizes $E$ as a Hamiltonian vector field, i.e.\ if $E=H_\Phi$. 
We should think of $\Phi$ as defining a matrix valued map
\[
\Phi_{ij}:\Rep(A,N)\r \Rep(k[t],N)=M_N(k)
\]
If $\Phi$ is a moment map for $\ldb -,-\rdb$ then $\Phi_{ij}$ is a moment
map for the induced Poisson bracket $\{-,-\}$. 
\subsection{Bisymplectic geometry}
\label{ref-2.4-10}
In the commutative case symplectic geometry is a special case of
Poisson geometry. A first version of non-commutative symplectic
geometry was introduced by Kontsevich in \cite{Kosymp}. See also
\cite{LebBock,Ginzburg1}. A version of Poisson geometry following
a similar philosophy was introduced by Crawley-Boevey in \cite{CB4}.

A different version of non-commutative symplectic geometry which
follows a similar philosophy as the Poisson geometry outlined in the
previous sections was introduced in \cite{CBEG} and baptized
``bisymplectic geometry''. We recall the definition. If $\delta\in D_A$
then we may define a double derivation
\[
i_\delta:\Omega A\r \Omega A \otimes \Omega A
\]
in the usual way. For $a\in A$ put
\[
i_\delta(a)=0\qquad i_\delta(da)=\delta(a)
\]
If $C$ is a graded $k$-algebra and $c=c_1\otimes \cdots \otimes c_n$ then we put
\begin{equation}
\label{ref-2.5-11}
{}^\circ c=\tau_{(1\cdots n)}(c)=(-1)^{|c_n|(|c_1|+\cdots+|c_{n-1}|)}c_nc_1\cdots c_{n-1}
\end{equation}
and if  $\phi:C\r C^{\otimes 2}$ is a linear map then we define
\[
{}^\circ \phi:C/[C,C]\r C: c\mapsto {}^\circ(\phi(c))
\]
We apply this with  $C=\Omega A$. 
Following \cite{CBEG} we put
\[
\imath_\delta={}^\circ\! i_\delta
\]
Also following \cite{CBEG} we say that an element $\omega\in \DR^2(A)$ is a \emph{bi-non-degenerate} if 
the map of $A$-bimodules
\[
\imath(\omega):D_{A}\r \Omega_{A}:\delta\mapsto \imath_\delta \omega
\]
is an isomorphism (in fact it is sufficient to assume surjectivity,
see Corollary \ref{ref-3.1.3-14} below). If in addition $\omega$ is closed in
$\DR(A)$ then we say that $\omega$ is \emph{bisymplectic}.

It was shown in \cite{CBEG} that if $\omega$ is bisymplectic the
$\tr(\omega)$ defines a symplectic form on representation spaces.

\section{Generalities about bimodules}
\label{ref-3-12}
In the framework of this paper the non-commutative version of a coherent sheaf
is a bimodule. In this short section we discuss some elementary
aspects of bimodules.
\subsection{Pairings}
Let $A$ be an arbitrary $k$-algebra. A pairing (or bilinear map) between
$A-A$ bimodules $P,Q$ is a map
\[
\langle-,-\rangle:P\times Q\r A\otimes A
\]
such that $\langle p,-\rangle$ is linear for the outer bimodule
structure on $A\otimes A$ and $\langle -,q\rangle$ is linear for the
inner bimodule structure on $A\otimes A$. The obvious example is of
course $P=Q^\ast$ and $\langle-,-\rangle$ is the evaluation pairing.
We say  that the pairing is non-degenerate if $P$, $Q$ are finitely
generated projective bimodules and the pairing induces an isomorphism
$Q\cong P^\ast$.

If $\langle -,-\rangle$ is a pairing between $P$ and $Q$ then the opposite
pairing between $Q$ and $P$ is given by 
\[
\langle q,p\rangle^\circ=\langle p,q\rangle''\otimes \langle p,q\rangle'
\]
If we have pairings between $P$ and $Q$ and between $P'$ and $Q'$ the
morphisms $\alpha:P\r P'$, $\beta:Q'\r Q$ are said to be \emph{adjoint} if
\[
\langle \alpha(p),q'\rangle=\langle p,\beta(q') \rangle
\]
We say that $\alpha$ is left adjoint to $\beta$ and that $\beta$ is a
right adjoint of $\alpha$. Note that if we change the pairings into
the opposite ones then a left adjoint becomes a right adjoint and vice
versa. Therefore we usually drop the left/right adjectives if what is
meant is clear from the context. If the pairing are non-degenerate
then we denote an adjoint often by $(-)^\ast$.

If we have a morphism $\alpha:P\r Q$ then then we say that $\alpha$
is anti-symmetric (or anti-self adjoint) if $-\alpha$ is left
adjoint to $\alpha$ for the appropriate pairings, i.e.\ if
 \[
\langle p,\alpha(p')\rangle=-\langle \alpha(p),p'\rangle^\circ
 \]
 for $p,p'\in P$. 

 Or written out explicitly
\begin{equation}
\label{ref-3.1-13}
\langle p,\alpha (p')\rangle'\otimes \langle p,\alpha (p')\rangle''
=-\langle p',\alpha (p)\rangle''\otimes \langle p',\alpha (p)\rangle'
\end{equation}

If we have a pairing between $P$ and $Q$ as above then $p\in P$ defines
a double derivation of degree $-1$
\[
i_p:T_A(Q)\r T_A(Q)\otimes T_A(Q)
\]
such that $i_p(q)=\langle p,q\rangle$.  Since $Q$ and $P$ are related
by the opposite pairing we may define $i_q$ for $q\in Q$ in the same way. I.e.
\[
i_q:T_A(P)\r T_A(P)\otimes T_A(P)
\]
is the double derivation of degree $-1$ such that $i_q(p)=\langle
q,p\rangle^\circ$.

We define $\imath_p={}^\circ i_p$ and $\imath_q={}^\circ i_q$ 
as in \S\ref{ref-2.4-10}. If $\omega\in T^2_A Q$ then
we put
\[
\imath(\omega):P\r Q:p\mapsto \imath_p(\omega)
\]
This is a bimodule morphism between $P$ and $Q$. 
The following property of $\imath(\omega)$ will be used below. 
\begin{propositions} The bimodule morphism  $\imath(\omega)$ is anti-symmetric.
\end{propositions}
\begin{proof} We need to prove that $\langle
  p,\imath(\omega)(p')\rangle=\langle p,\imath_{p'}\omega\rangle$
  satisfies \eqref{ref-3.1-13}.  We write $\omega$ as $\omega'\otimes
  \omega''$ with $\omega',\omega''\in Q$. Then we have
\[
i_{p'}(\omega)=i_{p'}(\omega')\omega''-\omega' i_{p'}(\omega'')
\]
so that we get
\begin{align*}
\imath_{p'}(\omega)&=i_{p'}(\omega')''\omega''i_{p'}(\omega')'
-i_{p'}(\omega'')''\omega'i_{p'}(\omega'')'\\
&=\langle p',\omega'\rangle''\omega''\langle p',\omega'\rangle'
-\langle p',\omega''\rangle''\omega'\langle p',\omega''\rangle'
\end{align*}
and hence
\begin{align*}
  \langle p,\imath_{p'}\omega\rangle&=\langle
  p',\omega'\rangle''\langle p,\omega''\rangle\langle p',\omega'\rangle'
  -\langle p',\omega''\rangle''\langle p,\omega'\rangle\langle p',\omega''\rangle'
\\
&=\langle
  p',\omega' \rangle''\langle p,\omega''\rangle'\otimes \langle p,\omega''\rangle''
\langle p',\omega'\rangle'
-\langle p',\omega''\rangle''\langle p,\omega'\rangle'\otimes 
\langle p,\omega'  \rangle''
\langle p',\omega''\rangle'
\end{align*}
which is clearly satisfies \eqref{ref-3.1-13}.
\end{proof}
\begin{corollarys} Assume that the pairing $\langle-,-\rangle$ is
non-degenerate. If $\imath(\omega)$ is surjective then it is an isomorphism.
\end{corollarys}
\begin{proof} Let $\alpha=\imath(\omega)$. If $\alpha$ is surjective
then its left adjoint $\alpha^\ast:Q^\ast\r P^\ast$ is injective. Since
the pairing is non-degenerate we have $Q^\ast=P$ and $P^\ast=Q$. Hence
we may identify $\alpha^\ast$ with $-\alpha$.
So $\alpha$ is both injective and surjective and hence it is an isomorphism.
\end{proof}
We can now state the following corollary which was asserted in \S\ref{ref-2.4-10}.
\begin{corollarys} \label{ref-3.1.3-14} Assume that  $\omega\in \DR^2(A)$ is such that the
map of bimodules 
\[
\imath(\omega):D_{A}\r \Omega_{A}:\delta\mapsto \imath_\delta \omega
\]
is surjective. Then $\omega$ is bi-non-degenerate. 
\end{corollarys}
If $\langle-,-\rangle:P\times Q\r A\otimes A$ then we call a \emph{dual
bases} sets of elements $p_\alpha\in P$, $q_\alpha\in Q$ such that
$\langle p_\alpha,-\rangle' q_\alpha \langle p_\alpha,-\rangle''$ is
the identity map on $Q$. It is easy to see that this equivalent to
$\langle -,q_\alpha\rangle'' p_\alpha \langle -,q_\alpha \rangle'$ being
the identity map on $P$. 
\begin{propositions} 
\label{ref-3.1.4-15} Assume that the pairing
  $\langle-,-\rangle:P\times Q\r A\otimes A$ is non-degenerate.  Then the
map
\[
\imath:T^2 _AQ\r \Hom_{A^e}(P,Q)
\]
defines an isomorphism between $(T_AQ/[T_AQ,T_AQ])_2$ and the anti-symmetric elements
of $\Hom_{A^e}(P,Q)$.
\end{propositions}
\begin{proof} It is easy to write down an explicit inverse to
  $\imath$. Let $p_\alpha\in P$, $q_\alpha\in Q$  be dual bases. 
  Then the element of $T^2_AQ$ corresponding to $\alpha\in
  \Hom_{A^e}(P,Q)^{\text{asym}}$ is given by
  $-\frac{1}{2}\alpha(p_\alpha)q_\alpha$.
\end{proof}
We will also have occasion to use the following result
\begin{propositions} \label{ref-3.1.5-16} Assume $\langle-,-\rangle$ is a
  non-degenerate pairing between $P$ and $Q$.
\begin{enumerate}
\item Let $\omega\in (T_AQ/[T_AQ,T_AQ])_n$. If $\imath_p(\omega)=0$ for all
  $p\in P$ then $\omega=0$.
\item Let $\eta\in T^nQ$. If the projection of $i_p(\eta)$
on $A\otimes T^{n-2}Q$ is zero for al $p\in P$ then $\eta=0$.
\end{enumerate}
\end{propositions}
\begin{proof} We select dual bases $p_\alpha$, $q_\alpha$.
The following formulas are easily verified for $\omega\in T^n_AQ$
\[
q_\alpha \imath_{p_\alpha}(\omega)=n\omega\qquad \mod [-,-]
\]
\[
(\pr_1(i_{p_\alpha}(\eta)))' q_\alpha (\pr_1(i_{p_\alpha}(\eta)))''=\eta
\]
From this the stated results follow.
\end{proof}
\subsection{Double Lie algebroids}
It will be convenient to make the following definition.
\begin{definitions}
\label{ref-3.2.1-17}
A \emph{double Lie algebroid} over $A$ is an $A$-bimodule $L$ together with
a (graded) double Poisson bracket of degree $-1$ on $T_AL$. 
\end{definitions}
The archetypical example of a double Lie algebroid is $D_A$ where we equip
$T_AD_A=DA$ with its double Schouten bracket (see \S\ref{ref-2.3-3}). 

Assume that $L$ is a double Lie algebroid with double bracket
$\ldb-,-\rdb$.  As for $DA$ we have associated operations
$\ldb-,-\rdb_l$, $\ldb-,-\rdb_r$ which are homomorphisms $L\times L\r L\otimes A$ and
$L\times L\r A\otimes L$ respectively that are defined by
\[
\ldb l_1,l_2\rdb=\ldb l_1,l_2\rdb_l+\ldb l_1,l_2\rdb_r
\]
$\ldb-,-\rdb_l$ and $\ldb-,-\rdb_r$ determine each other via
\[
\ldb l_1,l_2\rdb_r=-\ldb l_2,l_1\rdb_l^\circ
\]
It follows that the minimal data necessary to specify a double Lie
algebroid is given by
\begin{align*}
\ldb -,-\rdb_l&:L\times L\r L\otimes A\\
\ldb -,-\rdb&:L\times A\r A
\end{align*}
One can write down a minimal set of axioms these operations have to
satisfy (see \cite[(3.4-1)-(3.8-1)]{VdB33}). However it is often more
straightforward to use Definition \ref{ref-3.2.1-17} directly. 
\subsection{Representation spaces}
Let's now discuss how this plays out with representation spaces. 
if $P$ is  an $A$-bimodule then we define $P_N$ as the $A_N$ (cfr \S\ref{ref-2.1-0})-module
generated by symbols $p_{ij}$ which are linear in $p\in P$ and which
satisfy for $a\in A$:
\begin{align*}
(ap)_{ij}&=a_{iu}p_{uj}\\
(pa)_{ij}&=a_{uj}p_{iu}
\end{align*}
In this way we obtain an additive functor
\[
(-)_N:\Mod(A^e)\r \Mod(A_N)
\]
which sends finitely generated bimodules to finitely generated modules. 

This functor has a more intrinsic description as follows.
\begin{lemmas} \label{ref-3.3.1-18} \cite{CBEG} Consider $M_N(A_N)$ as an $A$-bimodule via the
$k$-algebra morphism $A\r M_N(A_N):a\mapsto:a_{ij}$. Consider $M_N(A_N)$
in addition as an $A_N$-module via the diagonal embedding $A_N\r M_N(A_N)$.
Then there is a natural isomorphism
\[
P_N\cong P\otimes_{A^e} M_N(A_N)
\]
In particular $(-)_N$ is right exact and sends projective bimodules to
projective modules.
\end{lemmas}
\begin{proof} We will content ourselves by giving the maps.
Let $f_{ij}$ be the usual elementary matrix in $M_N(A_N)$. Then one
isomorphism is given by
\[
P_N\r P\otimes_{A^e} M_N(A_N):p_{ij}\r p\otimes f_{ji}
\]
with the obvious definition for the inverse isomorphism.
\end{proof}
\begin{lemmas} \label{ref-3.3.2-19} We have $(T_A P)_N=S_{A_N} P_N$.
\end{lemmas}
\begin{proof} It is easy to see that both sides have the same generators
and relations. 
\end{proof}
\begin{lemmas} If $\langle-,-\rangle:P\times Q\r A\otimes A$ is
  non-degenerate then the corresponding pairing between $P_N$ and
  $Q_N$ (obtained by applying the standard index convention from
  \S\ref{ref-2.3-3})
\[
\langle-,-\rangle:P_N\otimes Q_N\r A: (p_{ij},q_{uv})\mapsto 
\langle p,q\rangle'_{uj}\langle p,q\rangle''_{iv}
\]
is non-degenerate as well.
\end{lemmas}
\begin{proof}
It follows from non-degeneracy that we may select 
$p_\alpha\in P$, $q_\alpha\in Q$ such that 
\[
q=\langle p_\alpha,q\rangle' q_\alpha \langle p_\alpha,q\rangle''
\]
for all $q\in Q$. 
It follows
\begin{align*}
\langle p_{\alpha,ij}, q_{uv} \rangle q_{\alpha,ji}
&=\langle p_{\alpha},q\rangle'_{uj}\langle p,q_\alpha\rangle''_{iv}  q_{\alpha,ji}\\
&=(\langle p_\alpha,q\rangle' q_\alpha \langle p_\alpha,q\rangle'')_{uv}\\
&=q_{uv}
\end{align*}
and hence the induced map 
$
Q_N\r P^\ast_N:q_{uv}\mapsto \langle -,q_{uv}\rangle
$
is injective  and split by $\phi\mapsto q_{\alpha,ji}\phi(p_{\alpha,ij})$.
Employing the dual argument we find that it is an isomorphism.
\end{proof}
\begin{propositions} 
\label{ref-3.3.4-20}
Assume that $A$ is smooth. Then the
  map $D_A\r \Der(A)$ defined by \eqref{ref-2.2-4} yields an isomorphism
  $(D_A)_N \r \Der(A)$.  In particular the tangent space to
  $\Rep(A,N)$ is generated by the vector fields $\delta_{ij}$ for
  $\delta\in \Der(A,N)$.
\end{propositions} 
\begin{proof}
It is easy to see that $(\Omega_A)_N=\Omega_{A_N}$.  By the previous lemma
we find 
\[
\Der(A_N)\cong\Omega_{A_N}^\ast\cong((\Omega_A)_N)^\ast\cong(D_A)_N
\]
It is easy to check that that the actual isomorphism is the asserted one. 
\end{proof}
\begin{corollarys} \label{ref-3.3.5-21} We have $(\Omega A)_N=\bigwedge_{A_N}
  \Omega_{A_N}$ (where $\Omega_{A_N}$ refers to the ordinary
  commutative differentials) and if $A$ is smooth then
  $(DA)_N=\bigwedge^\ast_{A_N} \Der(A_N)$.
\end{corollarys}
\begin{proof}
The statement about $\Omega A$ is easy. The statement about $DA$ follows from
\eqref{ref-3.3.4-20} and \eqref{ref-3.3.2-19}. 
\end{proof}
Below we sometimes use the map
\[
X:P\r M_N(P_N):p\mapsto (p_{ij})_{ij}
\]
For later use we mention some additional result. 
\begin{lemmas} \label{ref-3.3.6-22}  Let $P,Q,\langle-,-\rangle$ be as above and let $p\in P$.
Then the standard index convention (see \S\ref{ref-2.3-3}) applies to the operator $\imath_p$. I.e.
if $\omega\in T^2_A Q$ then we have
\[
i_{p_{ij}}(\omega_{uv})=i_p(\omega)'_{uj}i_p(\omega)''_{iv}
\]
In particular we obtain
$
i_{p_{ij}}(\tr(\omega))=(\imath_p \omega)_{ij}
$ or in more suggestive notation
\begin{equation}
\label{ref-3.2-23}
\tr(\omega)(X(p))=X(\imath(\omega)(p))
\end{equation}
where we view $\tr(\omega)$ as map from $P_N$ to $Q_N$. 
\end{lemmas}
\begin{lemmas} If $L$ is a double Lie algebroid over $A$ then $L_N$ is
a Lie algebroid over $A_N$ with Lie bracket and anchor map
given by the standard index convention. I.e.\ for $l,l_1,l_2\in L$, $a\in A$
\begin{align*}
[l_{1,ij},l_{2,uv}]&=\ldb l_1,l_2\rdb'_{uj}\ldb l_1,l_2\rdb''_{iv}\\
\rho(l_{ij})(a_{uv})&=\ldb l,a\rdb'_{uj}\ldb l,a\rdb''_{iv}
\end{align*}
\end{lemmas}
\section{Quasi-Poisson and quasi-Hamiltonian  $G$-spaces}
\label{ref-4-24}
Here we summarize some definitions from \cite{AKM,AMM}. 
Let $G$ be a linear algebraic group  and put $\frak{g}=\Lie(G)$. We
assume that $\frak{g}$ carries $G$-invariant symmetric bilinear form $(-,-)$. 

Let $(f_a)_a$,
$(f^a)_a$ be dual bases of $\frak{g}$. Then there is a canonical
invariant element $\phi\in \wedge^3\frak{g}$ given by
\[
\phi=\frac{1}{12} c^{abc} f_a\wedge f_b \wedge f_c
\]
where 
\[
c^{abc}=(f^a,[f^b,f^c])
\]
If $\xi \in \frak{g}$ we have left and right invariant vector fields
on $G$ defined by
\begin{align*}
(\xi^L(f))(x)&=\frac{d\ }{dt} f(x e^{\xi t})\\
(\xi^R(f))(x)&=\frac{d\ }{dt} f(e^{\xi t}x)
\end{align*}

Assume now that $G$ acts on a smooth affine variety $X$. If $\xi\in \frak{g}$
then $(\xi_x)_{x\in X}$ is the vector field on $X$ defined by 
\[
(v_{\xi_x}(f))(x)=\frac{d\ }{dt} f(e^{-t\xi} x)
\]
This convection is such that $\frak{g}\r TX:\xi\mapsto (\xi_x)_x$ is a morphism of Lie algebras.

The  element $\phi\in \wedge^3\frak{g}$ induces a three vector field
$\phi_X$ on $X$. Following \cite{AKM} an element $P\in \bigwedge^2_{\Oscr(X)} 
\Der(\Oscr(X))$ is said to be a \emph{quasi-Poisson bracket} if 
\[
\{P,P\}=\phi_X
\]
A \emph{Hamiltonian quasi-Poisson $G$-space} is a triple $(X,P,\Phi)$ 
such that $(X,P)$ is quasi-Poisson and such that $\Phi$ is a so-called
\emph{multiplicative moment map}. I.e. \
\[
\{h\circ \Phi,-\}=\frac{1}{2} f^a_X\left((f_a^L+f_a^R)(h)\circ \Phi\right)
\]
for all $f\in \Oscr(G)$.

A Hamiltonian quasi-Poisson $G$-space is said to be \emph{non-degenerate} 
if for $x\in X$ the map
\[
T^\ast_x\oplus \frak{g}\mapsto T_x:(\eta,\xi)\mapsto P_x(\eta)+\xi_x
\]
is surjective. 

\medskip

Now we let  $\theta$ and $\bar{\theta}$ be respectively
the left and right invariant $\frak{g}$-valued Maurer-Cartan forms:
\[
\theta=g^{-1}dg\qquad \bar{\theta}=dg\cdot g^{-1}
\]
We let $\chi$ be the canonical $G$-invariant three form on $G$:
\[
\chi=\frac{1}{12}(\theta,[\theta,\theta])
\]
A quasi-Hamiltonian $G$-variety is a triple $(X,\omega,\Phi)$ where
$M$ is a smooth $G$-variety, $\omega\in (T^{\ast,2})^G$ and $\Phi:X\r
G$ is a $G$-equivariant map (for the given action on $X$ and the adjoint
action of $G$ on $G$) such that the following axioms are satisfied:
\begin{enumerate}
\item[(B1)] $d\omega=\Phi^\ast \chi$.\footnote{We follow the
    convention from \cite[Def 10.1]{AKM}. In \cite{AMM} this formula is given
    with a minus sign}
\item[(B2)] $\forall \xi\in \frak{g}:
  i_{\xi_X}(\omega)=\frac{1}{2}\Phi^\ast(\theta+\bar{\theta},\xi)$.
\item[(B3)] For all $x$ in $X$ we have
\[
\ker \omega_x=\{\xi_x\mid \xi \in \ker(\Ad_{\Phi(x)}+1)\subset \frak{g}\}
\]
\end{enumerate}
Let us explain how to read (B3). By definition $\omega_x$ is the map
\[
T_X\r T^\ast_X:\delta\mapsto i_\delta(\omega)
\]
evaluated in $x$.  For $\xi\in \frak{g}$, $\xi_x$
is the vector field on $X$ given by $\xi$ evaluated in $x$. (B3)
states that $\ker \omega_x$ is a specific part of the image
of the map $\frak{g}\r T_{X,x}:\xi\mapsto \xi_x$. 
As the form of condition (B3) is not so convenient for us so we state an 
equivalent version.
\begin{enumerate}
\item[(B3')]  The map
\[
T_x\oplus \frak{g}\r T^\ast_x:(\delta,\xi)\mapsto \omega_x(\delta)
+ (\xi, \Phi^\ast(\theta)_x)
\]
is surjective. 
\end{enumerate}
\begin{lemma} (B3) and (B3') are equivalent. 
\end{lemma}
\begin{proof} (Sketch) We first dualize (B3'). Let
  $\beta:T_x\r \frak{g}:\delta\mapsto i_{\delta}(\Phi^\ast(\theta))$. Then
(B3') is equivalent to the condition
\begin{enumerate}
\item[(B3'')] The map
$
(\omega_x,\beta):T_x\r T^\ast_x\oplus \frak{g} 
$
is injective. 
\end{enumerate}
It is a straightforward verification using (B2) that the following diagram
\begin{equation}
\label{ref-4.1-25}
\begin{CD}
 \frak{g} @>\alpha >> \frak{g}\\
@A\beta AA   @AA\gamma A\\
T_x @>>\omega_x> \Omega_x
\end{CD}
\end{equation}
with $ \alpha(\xi)=-\frac{1}{2}(1+\Ad(\Phi(x)))(\xi)$ and $(
\gamma(\eta),\xi)=i_{\xi_x}(\eta)$ is commutative, Furthermore using
the properties of $\theta$, $\bar{\theta}$ (see e.g\ \cite[Ch
II]{Helgason}) we find $\beta(\xi_x)=(1-\Ad(\Phi(x))^{-1})(\xi)$.

Let us show that (B3) implies (B3'').  Assume that $\delta$ is such that $i_\delta(\omega_x)=0$
and $i_{\delta}(\Phi^\ast(\theta))=0$. By (B3) we have that
$\delta=\xi_x$ where $(1+\Ad(\Phi(x)))(\xi)=0$. Since $\beta(\xi_x)=0$
we also have $(1-\Ad(\Phi(x)))(\xi)=0$ by diagram \eqref{ref-4.1-25}. These two facts together imply
that $\xi=0$. 

Now we prove the converse. Assume that $\delta\in \ker(\omega_x)$. Then $\xi=\beta(\delta)\in \ker \alpha$.
I.e.
\begin{equation}
\label{ref-4.2-26}
(1+\Ad(\Phi(x)))(\xi)=0
\end{equation}
Then 
\[
\beta(\xi_x)=(1-\Ad(\Phi(x)^{-1})(\xi)=2\xi
\]
Since $(1+\Ad(\Phi(x)))(\xi)=0$ we also have $\omega_x(\xi_x)=0$.

It follows that $=\delta-\frac{1}{2}\xi_x$ is both in the kernel of
$\beta$ an $\omega_x$. Hence by (B3'') $\delta=\frac{1}{2}{\xi_x}$.
(B3) now follows from \eqref{ref-4.2-26}.
\end{proof}

We now state the main theorem of \cite[\S 10]{AKM}. Write $\theta=f^a \theta_a$
where $\theta_a\in T_X^\ast$ and similarly for $\bar{\theta}$. 
\begin{theorem} \label{ref-4.2-27} \cite[Thm 10.3]{AKM} Every non-degenerate Hamiltonian
  quasi-Poisson space $(X,P,\Phi)$ carries a unique 2-form $\omega$
  such that $(X,\omega,\Phi)$ is a quasi-Hamiltonian $G$-space and such
that $\omega$ and $P$ satisfy the following compatibility condition
\begin{equation}
\label{ref-4.3-28}
P\circ \omega=1-\frac{1}{4} f^{a}_{X} \otimes \Phi^\ast(\theta_a-\bar{\theta}_a)
\end{equation}
(as maps $T^\ast_X\r T_X$). Conversely on every quasi-Hamiltonian $G$-space
$(X,\omega,\Phi)$ there is a unique bivector field $P$ such that $(M,X,\Phi)$
is a non-degenerate quasi-Hamiltonian $G$-manifold and \eqref{ref-4.3-28} 
is satisfied. 
\end{theorem}
We will say that $(P,\omega,\Phi)$ are compatible 
if \eqref{ref-4.3-28} holds.

\section{Hamiltonian double quasi-Poisson algebras}
The non-commutative version of quasi-Poisson algebras was worked out
in \cite{VdB33}. 
\begin{convention}
From now on our non-commutative algebras are always smooth.
\end{convention}
Let $P\in (DA/[DA,DA])_2$. We say that $P$ is a \emph{double quasi-Poisson bracket} if
the following condition holds
\begin{enumerate}
\item[($\mathbb{P}$1)] $\{P,P\}=\frac{1}{12} E^3 \qquad \mod [-,-]$.
\end{enumerate}
In addition
an invertible element $\Phi\in A$ is said to be a \emph{multiplicative
moment map} if the following condition holds.
\begin{enumerate}
\item[($\mathbb{P}$2)] $\{\Phi,-\}=\frac{1}{2}(\Phi E+E\Phi)$.
\end{enumerate}
Finally we say that $P$ is \emph{non-degenerate} if the following
condition holds
\begin{enumerate}
\item[($\mathbb{P}$3)] The map $\Omega_A\oplus AEA\r D_A :(\eta,\delta)\mapsto \imath(P)(\eta)+\delta$ is surjective.
\end{enumerate}
It was proved in \cite{VdB33} that ($\PP$1) and ($\PP$2) imply the
corresponding properties on representation spaces. The same is
also true for ($\PP$3) as the following proposition shows. 
\begin{proposition} Assume that ($\PP$1)($\PP$3) hold.
Then $\tr(P)$ defines a non-degenerate quasi-Poisson bracket on $\Rep(A,N)$
\end{proposition}
\begin{proof}
  Applying the right exact functor $(-)_N$ (lemma \ref{ref-3.3.1-18}) and using
  Proposition \ref{ref-2.3.1-9} together with Corollary \ref{ref-3.3.5-21}
  we find that the map
\[
\Omega_{A_N}\oplus A_N \otimes \frak{gl}_N\r \Der(A_N)
\]
is surjective. This finishes the proof. 
\end{proof}
The following result is a non-commutative version of \cite[Thm.\
2.5]{BuCr}. 
\begin{theorem} \label{ref-5.3-29} Assume that $(A,P)$ is a double quasi-Poisson algebra.
  Then $\tilde{\Omega}_A=\Omega_A\oplus AEA$ has the structure of a
  double Lie algebroid where the double bracket is defined as follows.
\begin{align*}
\ldb da,b\rdb_{\tilde{\Omega}A}&=\ldb a,b\rdb\\
\ldb da,db\rdb_{\tilde{\Omega}A}&= d\ldb a,b\rdb+
\frac{1}{4} [b,[a,E\otimes 1-1\otimes E]_\ast]\\
\ldb E,X\rdb_{\tilde{\Omega}A}&=X\otimes 1-1\otimes X
\end{align*}
for $a,b\in A$, $X\in T_A\tilde{\Omega}_A$ and where $[-,-]_\ast$
denotes the commutator for the inner $A$-bimodule structure on
$AEA\otimes AEA$. Furthermore the map
\[
\Omega_A\oplus AEA\r D_A 
\]
defined in ($\PP$3) is a morphism of double Lie algebroids. 
\end{theorem}
Our proof of this theorem is a rather painful direct computation. We
omit the details.  The fact that ($\PP$3) defines a morphism of double
Lie algebroids translates into the following proposition which is an
analogue of \cite[Prop.\ 3.5.1]{VdB33}.
\begin{proposition} \label{ref-5.4-30} The following are equivalent
\begin{enumerate}
\item $\ldb-,-\rdb$ is a double quasi-Poisson bracket.
\item The following identity holds for all $a,b\in A$:
\[
\ldb H_a,H_b\rdb_l -H_{\ldb a,b\rdb}'\otimes \ldb a,b\rdb''=
\frac{1}{4} [b,[a,E\otimes1]_\ast]
\]
\item The following identity holds for all $a,b\in A$:
\[
\ldb H_a,H_b\rdb_r -\ldb a,b\rdb'\otimes H_{\ldb a,b\rdb''}=
-\frac{1}{4} [b,[a,1\otimes E]_\ast]
\]
\item The following identity holds for all $a,b\in A$:
\[
\ldb H_a,H_b\rdb -H_{\ldb a,b\rdb}=
\frac{1}{4} [b,[a,E\otimes1-1\otimes E]_\ast]
\]
where we use the convention $H_{x'\otimes x''}=H_{x'}\otimes x''+x'\otimes H_{x''}$.
\end{enumerate}
\end{proposition}
\begin{proof}
Below it  will be convenient to use the notation ${}^\circ(-)$ introduced in
\eqref{ref-2.5-11} as well as the convention 
\[
(a_1\otimes \cdots \otimes a_m)(b_1\otimes \cdots \otimes b_n)=
a_1\otimes \cdots \otimes a_m b_1\otimes \cdots b_n
\]
and similarly for longer products.
According to the proof of \cite[Prop. 3.5.1]{VdB33} we have for $a,b,c\in A$
\[
\ldb a,b,c\rdb=\tau_{23}((\ldb H_a,H_b\rdb_l -H_{\ldb a,b\rdb'}(-)\otimes \ldb a,b\rdb'')(c))
\]
Also according to \cite[\S5]{VdB33} a bracket is quasi-double Poisson
if it satisfies
\[
\ldb -,-,-\rdb=\frac{1}{12}\ldb -,-,-\rdb_{E^3}
\]
By the formulas in \cite[\S4]{VdB33} it follows 
\[
\ldb -,-,-\rdb_{E^3}=3 {}^\circ (E(a)^\circ E(b)^\circ E(c)^\circ)
\]
Hence
\begin{equation}
  \label{ref-5.1-31}
  (\ldb H_a,H_b\rdb_l -H_{\ldb a,b\rdb'}(-)\otimes \ldb a,b\rdb'')(c)
  =\frac{1}{4} \tau_{23} {}^\circ (E(a)^\circ E(b)^\circ E(c)^\circ)
\end{equation}
A straightforward verification yields
\[
(E(a)^\circ E(b)^\circ E(-)^\circ)=[b,[a,E\otimes 1]_\ast](c)
\]
This proves the equivence between (1) and (2). The equivalence between
(2) and (3) follows easily from the fact that $\ldb a,b\rdb_r=-\ldb
b,a\rdb^\circ _l$. The sum of (2) and (3) yields (4). To go back we use
projection. 
\end{proof}
\section{Quasi-bisymplectic algebras}
\label{ref-6-32}
The algebras we will introduce are a non-commutative analogue of
quasi-Hamiltonian $G$-spaces introduced in \cite{AMM} (see \S\ref{ref-4-24}).
By definition a \emph{quasi-bisymplectic algebra} will be a triple
$(A,\omega,\Phi)$ where $\omega\in \DR^2(A)$ and $\Phi\in A^\ast$
satisfying the following conditions
\begin{enumerate}
\item[($\mathbb{B}$1)] $d\omega=\frac{1}{6} (\Phi^{-1} d\Phi)^3\quad \mod [-,-]$.
\item[($\mathbb{B}$2)] $\imath_{E}\omega=\frac{1}{2} (\Phi^{-1}d\Phi+d\Phi\cdot \Phi^{-1})$
\item[($\mathbb{B}$3)] The  map
\[
D_A\oplus Ad\Phi A\r \Omega_A :(\delta,\eta)\mapsto
\imath(\omega)(\delta)+\eta
\]
is surjective.
\end{enumerate}
\begin{proposition} If $(A,\omega,\Phi)$ is a quasi-bisymplectic
algebra then $(\Rep(A,N),\tr(\omega),X(\Phi))$ is a quasi-Hamiltonian $\Gl_N$-space. 
\end{proposition}
\begin{proof} It is easy to see that we have
\begin{equation}
\label{ref-6.1-33}
 X(\Phi^{-1}d\Phi)=X(\Phi)^\ast(\theta)\qquad  X(d\Phi\cdot \Phi^{-1})=X(\Phi)^\ast(\bar{\theta})
\end{equation}
 and 
\begin{equation}
\label{ref-6.2-34}
X(\Phi)^\ast(\chi)=\frac{1}{6} \tr(\Phi^{-1}d\Phi)^3
\end{equation}
This implies  (B1). To check (B2) we test it with $\xi=f_{ji}$. Then
according to lemma \ref{ref-3.3.6-22} and Prop.\ \ref{ref-2.3.1-9} we have 
\[
i_{f_{ji}}(\tr(\omega))=(\imath_E \omega)_{ij}
\]
We also have 
\[
X(\Phi)^\ast(\theta+\bar{\theta},f_{ji})=(\Phi^{-1} d\Phi+d\Phi\cdot \Phi^{-1})_{ij}
\]
so that we see that ($\mathbb{B}$2) implies (B2). 

Now we consider (B3) (or rather its variant (B3')). Apply the functor $(-)_N$
to ($\mathbb{B}$3). Using  Corollary \ref{ref-3.3.5-21}
  we find that the map
\[
\Der(A_N)\oplus \sum_{ij} A_N (d\Phi)_{ij}\r \Omega_{A_N}
\]
is surjective. But the forms $(d\Phi)_{ij}$ generate the same 
$A_N$-module as
$X(\Phi)^\ast(\theta)_{ij}=(\Phi^{-1}d\Phi)_{ij}$, finishing
the proof.
\end{proof}
\begin{remark} Note that the appearance of the element
  $(\Phi^{-1}d\Phi)^3\in \DR^3(A)$ is rather natural in view of
  Example \ref{ref-2.1-2}. This is the only non-zero element in $H^3(\DR(A))$ that can be constructed from the single element
  $\Phi$.
\end{remark}
Let us say that $P\in D^2A$, $\omega\in \Omega^2 A$, $\Phi\in A^\ast$
are \emph{compatible} if the following identity holds for all
$\delta\in D_A$:
\begin{equation}
\label{ref-7.3-32} \tag{$\CC$}
(\imath(P)\circ
i(\omega))(\delta)=\delta-\frac{1}{4}\delta(\Phi)''(E\Phi^{-1}-\Phi^{-1}E)\delta(\Phi)'
\end{equation}
\begin{proposition}
If $P,\omega,\Phi$ are compatible then
  $(\tr(P),\tr(\omega),X(\Phi))$ are compatible as well (see \eqref{ref-4.3-28}).
\end{proposition}
\begin{proof}
We apply the map $X(-)$. Using Cor.\ \ref{ref-3.3.5-21}, \eqref{ref-3.2-23} and
Prop.\ \ref{ref-2.3.1-9} we find 
\begin{align*}
  (\tr(P)\circ \tr(\omega))(\delta_{ij})&=\delta_{ij} -
  \frac{1}{4}\delta(\Phi)''_{iu}(f_{vu} (\Phi^{-1})_{vw}-
  (\Phi^{-1})_{uv} f_{wv}) \delta(\Phi)'_{wj}\\
&=\delta_{ij} -
  \frac{1}{4}(f_{vu} (\Phi^{-1})_{vw}-
  (\Phi^{-1})_{uv} f_{wv}) \delta(\Phi)'_{wj}\delta(\Phi)''_{iu}\\
&=\delta_{ij} -
  \frac{1}{4}(f_{vu} (\Phi^{-1})_{vw}-
  (\Phi^{-1})_{uv} f_{wv}) \delta_{ij}(\Phi_{wu})\\
&=\delta_{ij} -
  \frac{1}{4}(f_{vu} (\Phi^{-1})_{vw} i_{\delta_{ij}}(d\Phi_{wu}) -
 i_{\delta_{ij}}(d\Phi_{wu}) (\Phi^{-1})_{uv} f_{wv}) \\
&=(\Id-\frac{1}{4} f_{vu}( (\Phi^{-1}d\Phi)_{vu} -
(d\Phi\cdot \Phi^{-1})_{vu} 
))(\delta_{ij})
\end{align*}
where we have viewed the $f_{ij}$ as vector fields on $\Rep(A,N)$.  We
are now done by \eqref{ref-6.1-33}.
\end{proof}
\section{Compatibility}
Our aim is to prove a non-commutative analogue of Theorem \ref{ref-4.2-27}.
\begin{theorem} \label{ref-7.1-35} Fix $\Phi\in A^\ast$. 
\begin{enumerate}
\item For every $P\in (DA/[DA,DA])_2$ satisfying ($\mathbb{P}$2)($\mathbb{P}$3)
there exists a unique $\omega\in (\Omega A/[\Omega A,\Omega A])_2$ 
satisfying  ($\mathbb{B}$2)($\mathbb{B}$3) and  \eqref{ref-7.3-32}.
\item For every $\omega\in (\Omega A/[\Omega A,\Omega A])_2$ satisfying ($\mathbb{B}$2)($\mathbb{B}$3)
there exists a unique $P\in (DA/[DA,DA])_2$ 
satisfying  ($\mathbb{P}$2)($\mathbb{P}$3) and  \eqref{ref-7.3-32}.
\item If $\omega$, $P$ correspond to one another as in (1)(2)
  then the integrability conditions ($\mathbb{B}$1) and ($\mathbb{P}$1) are
equivalent.
\end{enumerate}
\end{theorem}
We will prove this theorem in this section.  Throughout we  fix $\Phi\in A^\ast$ and we let
$\omega\in (\Omega A/[\Omega A,\Omega A]))_2$, $P\in (D A/[D A, D A]))_2$
be such that ($\mathbb{B}$2) and ($\mathbb{P}$2) are satisfied. 
First we consider the
following diagram which  summarizes a number of relevant maps 
which we will use below. 
\begin{equation}
\label{ref-7.1-36}
\xymatrix{
\Omega_A \ar[r]^{e}\ar[d]_{\imath(P)}
& AE^\ast A  \ar[d]^{\jmath} \ar[r]^{T^0} &Ad\Phi A\ar[r]^{c}
\ar[d]^{\imath}& 
\Omega_A\ar[d]^{\imath(P)}\\
D_A\ar[r]_{e} \ar[d]_{\imath(\omega)} & A(d\Phi)^\ast A \ar[d]^{\jmath}\ar[r]_{S^0}& AEA\ar[d]^{\imath} \ar[r]_{c} & D_A\ar[d]^{\imath(\omega)} \\
\Omega_A \ar[r]_{e}
& AE^\ast A \ar[r]_{T^0} &
Ad\Phi A\ar[r]_{c} 
& 
\Omega_A
}
\end{equation}
The following conventions are used: first of all we view $AE^\ast A$,
$A(d\Phi)^\ast A$, $AEA$ and $A d\Phi A$ as free bimodules 
with one generator. Furthermore we regard the diagram as being doubly
infinite with period $(3,2)$.  The bimodules occur in pairs related by
an obvious non-degenerate pairing: $(D_A,\Omega_A)$, $(AE^\ast A,
AEA)$ and $(A(d\Phi)^\ast A, A(d\Phi) A)$. 
 We now define the maps.
\begin{enumerate}
\item $c$ stands for  ``canonical map''.
\item $e$ is adjoint to $c$. The formulas for the two variants are as follows.
\begin{align*}
e(d\phi)&=\phi E^\ast-E^\ast \phi\\
e(\delta)&=\delta(\Phi)''(d\Phi)^\ast \delta(\Phi)'
\end{align*}
\item $\imath$ is the restricted version of $\imath(P)$ and
  $\imath(\omega)$. The two variants are given by ($\mathbb{P}$2)
  and ($\mathbb{B}$2). I.e.\ explicitly
\begin{align*}
\imath(d\Phi)&=\frac{1}{2}(E\Phi+\Phi E)\\
 \imath(E)&=\frac{1}{2}(\Phi^{-1}d\Phi+d\Phi\cdot \Phi^{-1})
\end{align*}
\item $\jmath$ is adjoint to
$-\imath$. The two variants  are given by the 
following formulas
\begin{align*}
\jmath(d\Phi^\ast)&=-\frac{1}{2}(\Phi^{-1}E^\ast+E^\ast \Phi^{-1})\\
\jmath(E^\ast)&=-\frac{1}{2}(\Phi
(d\Phi)^\ast+(d\Phi)^\ast \Phi)
\end{align*}
\item The maps $S^0$ and $T^0$ are adjoint to one another. They are
respectively given by the following formulas. 
\begin{align*}
S^0((d\Phi)^\ast) &=E\Phi^{-1}-\Phi^{-1} E\\
T^0(E^\ast )&=\Phi^{-1}d\Phi-
d\Phi \Phi^{-1}
\end{align*}
\end{enumerate}
From the stated adjointness properties it follows that the 
\eqref{ref-7.1-36} is  self dual, up to sign.  We have the following fact.
\begin{lemma} 
\label{ref-7.2-37}
\begin{enumerate}
\item The diagram \eqref{ref-7.1-36} is commutative.
\item Consider the diagram as doubly infinite. For any configuration of arrows 
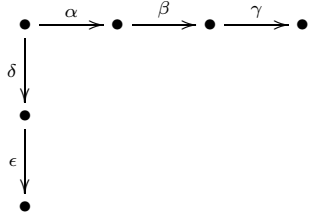
\[
\xymatrix{%
\bullet\ar[r]^\alpha\ar[d]_\delta &\bullet\ar[r]^\beta &\bullet\ar[r]^\gamma &\bullet\\
\bullet\ar[d]_\epsilon\\
\bullet}
\]
such that the vertical arrows do not involve $\imath(\omega),\imath(P)$ we
have 
\begin{equation}
\label{ref-7.2-38}
\frac{1}{4}\gamma \beta\alpha+ \epsilon\delta=1
\end{equation}
\item If we contract in \eqref{ref-7.1-36} the pairs of horizontal consecutive arrows that
  have $\Omega_A$ or $D_A$ in the middle then the $2\times 2$ subdiagrams
  in the resulting diagram are bicartesian.
\item If ($\mathbb{B}$3) holds then any subdiagram 
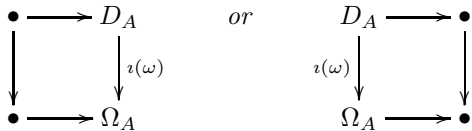
\[
\xymatrix{%
\bullet \ar[d]\ar[r] &D_A\ar[d]^{\imath(\omega)}\\
\bullet \ar[r] &  \Omega_A
}
\qquad \text{or}\qquad
\xymatrix{%
D_A \ar[d]_{\imath(\omega)}\ar[r] &\bullet\ar[d]\\
\Omega_A \ar[r] &  \bullet
}
\]
of \eqref{ref-7.1-36} is bicartesian.  
\item Similarly if ($\mathbb{P}$3) holds then any subdiagram 
\[
\xymatrix{%
\bullet \ar[d]\ar[r] &\Omega_A\ar[d]^{\imath(P)}\\
\bullet \ar[r] &  D_A
}\qquad\text{or}\qquad
\xymatrix{%
\Omega_A \ar[d]_{\imath(P)}\ar[r] &\bullet\ar[d]\\
D_A \ar[r] &  \bullet
}
\]
of \eqref{ref-7.1-36} is bicartesian. 
\end{enumerate}
\end{lemma}
\begin{remark} To avoid confusion: statements in the above lemma
which do not refer to $\imath(\omega)$ or $\imath(P)$ remain true if
these maps are not present in the diagram.
\end{remark}
\begin{proof} We leave the verification of (1) and (2)
 to the reader.  For (3) will only give an example.  E.g.\ we need to check in
particular that the following diagram is bicartesian
\[
\begin{CD}
Ad\Phi A @>ec>> AE^\ast A\\
@V\imath VV @VV \jmath V\\
AEA @>>ec> A(d\Phi)^\ast A
\end{CD}
\]
Put $R=A^e$, $\Phi_1=\Phi\otimes 1$, $\Phi_2=1\otimes \Phi$. The latter are
two commuting invertible elements of $R$. Removing all
confusing decoration the above diagram basically looks like
\[
\begin{CD}
R @>\Phi_1-\Phi_2>> R\\
@V \Phi_1+\Phi_2 VV @VV \Phi_1+\Phi_2  V\\
R @>>\Phi_1-\Phi_2 > R
\end{CD}
\]
It is now clear that this diagram is bicartesian (using the fact that
$\Phi_1$ and $\Phi_2$ are invertible). The other $2\times 2$-diagrams
in \eqref{ref-7.1-36} are similar.

(4) and (5) are similar. To prove (4) we only need to consider the
first diagram as the second follows by duality. Assume that there are
$\delta\in D_A$, $\eta\in Ad\Phi A$ such that
\begin{equation}
\label{ref-7.3-39}
\imath(\omega)(\delta)=c(\eta)
\end{equation} 
 We need
to prove that there is a unique $\Delta\in AEA$ such that $c(\Delta)=\delta$,
$\imath(\Delta)=\eta$. We need to worry only about existence since
(2) implies that 
\[
(ec,\imath):AEA\r A(d\Phi)^\ast A\oplus Ad\Phi A 
\]
is injective and hence $(c,\imath)$ is also injective.

Applying $e$ to \eqref{ref-7.3-39} we find $\jmath
e(\delta)=e\imath(\omega)(c)=ec(\eta)$. By (3) there is an element
$\delta'\in AEA$ such that $\imath(\delta')=\eta$,
$ec(\delta')=e(\delta)$. Replacing $\eta$ by $\eta-\imath(\delta')$,
$\delta$ by $\delta-c(\delta')$ we reduce to the case
$\eta=0$, $e(\delta)=0$, $\imath(\omega)(\delta)=0$.

Now we note that ($\mathbb{B}$3) implies by duality that
the following map is injective
\[
(\imath(\omega), e):D_A\r \Omega_A\oplus A(d\Phi)^\ast A
\]
Hence $\delta=0$ and we can take $\Delta=0$. 
\end{proof}
Using the maps in \eqref{ref-7.1-36} the compatibility condition
\eqref{ref-7.3-32} can be reformulated as 
\begin{equation}
\label{ref-7.4-40}
\imath(P)\imath(\omega)=1-\frac{1}{4}S
\end{equation}
where $S=cS^0e$. In fact it will be convenient to give some other equivalent
formulations of \eqref{ref-7.4-40}. Consider the following matrices
\[
\bar{\omega}=\begin{pmatrix}
\imath(\omega)& c\\
\frac{1}{4} S^0e & -\imath
\end{pmatrix}
\qquad
\bar{P}=\begin{pmatrix}
\imath(P)& c\\
\frac{1}{4}T^0e& -\imath
\end{pmatrix}
\]
and view them as maps between $D_A\oplus Ad\Phi A$ and $\Omega_A\oplus AEA$.

Furthermore consider the matrices
\[
\tilde{\omega}=\begin{pmatrix}
\imath(\omega)& 1\\
\frac{1}{4} S & -\imath(P)
\end{pmatrix}
\qquad
\tilde{P}=\begin{pmatrix}
\imath(P)& 1\\
\frac{1}{4}T& -\imath(\omega)
\end{pmatrix}
\]
and view them as maps between $D_A\oplus \Omega_A$ and $\Omega_A\oplus D_A$.
\begin{proposition} \label{ref-7.4-41} The following conditions are all equivalent.
\begin{align}
\label{ref-7.5-42} \imath(P)\imath(\omega)&=1-\frac{1}{4}S\\
\label{ref-7.6-43}  \imath(\omega)\imath(P)&=1-\frac{1}{4}T\\
\label{ref-7.7-44} \bar{P}\bar{\omega}&=\Id\\
\label{ref-7.8-45} \bar{\omega}\bar{P}&=\Id\\
\label{ref-7.9-46} \tilde{P}\bar{\omega}&=\Id\\
\label{ref-7.10-47} \tilde{\omega}\bar{P}&=\Id
\end{align}
\end{proposition}
\begin{proof}
\eqref{ref-7.5-42} and \eqref{ref-7.6-43} are equivalent since they are adjoint. To
understand \eqref{ref-7.7-44} we compute the product
\[
\bar{P}\bar{\omega}=\begin{pmatrix}
\imath(P)& c\\
\frac{1}{4}T^0e& -\imath
\end{pmatrix}
\begin{pmatrix}
\imath(\omega)& c\\
\frac{1}{4} S^0 e & -\imath
\end{pmatrix}
=
\begin{pmatrix}
\imath(P)\imath(\omega)+\frac{1}{4}S&\imath(P)c-c \imath \\
\frac{1}{4}(T^0e\imath(\omega) -\imath S^0 e)& 
\frac{1}{4} T^0 ec +\imath\imath
\end{pmatrix}
=
\begin{pmatrix}
\imath(P)\imath(\omega)+\frac{1}{4}S&0\\
0& 
1
\end{pmatrix}
\]
where we have used lemma \ref{ref-7.2-37}. It is now clear
that \eqref{ref-7.5-42} and \eqref{ref-7.7-44} are equivalent.  Similarly
\eqref{ref-7.6-43} and \eqref{ref-7.8-45} are equivalent. 
To understand \eqref{ref-7.9-46}
we write the product out again explicitly.
\[
\tilde{P}\tilde{\omega}=
\begin{pmatrix}
\imath(P)& 1\\
\frac{1}{4}T& -\imath(\omega)
\end{pmatrix}
\begin{pmatrix}
\imath(\omega)& 1\\
\frac{1}{4} S & -\imath(P)
\end{pmatrix}
=
\begin{pmatrix}
\imath(P)\imath(\omega)+\frac{1}{4} S&0\\
0&\imath(\omega) \imath(P)+\frac{1}{4} T
\end{pmatrix} 
\]
where we have used lemma \ref{ref-7.2-37} again. It follows that
\eqref{ref-7.9-46} is equivalent to \eqref{ref-7.5-42}\eqref{ref-7.6-43} simultaneously. The same holds for \eqref{ref-7.9-46}.
\end{proof}
\begin{proof}[Proof of Theorem \ref{ref-7.1-35}(1)]

 If  $\imath(\omega)$ exists it
will have the following two properties for $\eta\in \Omega_A$, $\delta\in AEA$
\begin{align}
\label{ref-7.11-48} \imath(\omega)\imath(P)(\eta)&=\eta-\frac{1}{4}T(\eta)\\
\label{ref-7.12-49} \imath(\omega)c(\delta)&=c\imath(\delta)
\end{align}
where we have used lemmas \ref{ref-7.2-37} and Prop.\ \ref{ref-7.4-41}. 
Since we are assuming ($\mathbb{P}$3)  every element of $D_A$ can be
written as a sum $\imath(P)(\eta)+c(\delta)$ for $\eta\in \Omega_A$
and $\delta\in AEA$. Hence the properties \eqref{ref-7.11-48}\eqref{ref-7.12-49}
characterize $\imath(\omega)$
uniquely (and $\omega$ as well by  Prop.\
   \ref{ref-3.1.4-15}). 

We still need to prove two things. 
\begin{enumerate}
\item The equations \eqref{ref-7.11-48}\eqref{ref-7.12-49} are non-contradictory.
\item The resulting $\imath(\omega)$ is actually anti-symmetric (so that
it genuinely comes from an element $\omega\in \Omega^2_A$). 
\item $\imath(\omega)$ satisfies ($\mathbb{B}$3).
\end{enumerate}
We discuss (1) first. Suppose that $\imath(P)(\eta)=c(\delta)$. Then
by lemma \ref{ref-7.1-36} there exist $\eta_0\in Ad\Phi_pA$ such that
$\eta=c(\eta_0)$, $\delta=\imath(\eta_0)$. Then
\begin{align*}
\eta-\frac{1}{4}T(\eta)&=c(\eta_0)-\frac{1}{4}c T^0 e c(\eta_0)\\
&=c(\eta_0-\frac{1}{4} T^0 e c(\eta_0))\\
&=c\imath\imath(\eta_0)\\
&=c\imath(\delta)
\end{align*}
where we have used lemma \ref{ref-7.1-36}. 

Now we prove (2). To avoid confusion we write $X$ for the map
$\imath(\omega)$ we have constructed.  Thus we have
\begin{equation}
\label{ref-7.13-50}
\begin{split}
X\imath(P)&=1-\frac{1}{4}T\\
Xc&=c\imath
\end{split}
\end{equation}
Dualizing these equations we get
\begin{equation}
\label{ref-7.14-51}
\begin{split}
-\imath(P) X^\ast&=1-\frac{1}{4}S\\
eX^\ast &=-\jmath e
\end{split}
\end{equation}
We write \eqref{ref-7.13-50} in matrix form.
\[
X(\imath(P)\, c)=(1-\frac{1}{4}T\, c\imath)
\]
Left multiplication by $(i(P)\, -e)^t$ yields
\begin{align*}
\begin{pmatrix} 
i(P)\\
-e
\end{pmatrix}
X(\imath(P)\, c)&=\begin{pmatrix} \imath(P)\\\, -e\end{pmatrix}
(1-\frac{1}{4} T\, c\imath)\\
&=\begin{pmatrix} 1-\frac{1}{4}S \\ -\jmath e\end{pmatrix}
(\imath(P)\, c)
\end{align*}
Since $(\imath(P)\, c)$ is injective we conclude 
\begin{align*}
\imath(P)(X)&=1-\frac{1}{4}S \\
-eX&=-\jmath e
\end{align*}
Comparing with \eqref{ref-7.14-51}, and using that $(\imath(P) -e)$ is
surjective we conclude that indeed $X^\ast=-X$.

We finish the proof by noting that (3) follows immediately from
\eqref{ref-7.6-43}.
\end{proof}
\begin{remark} Assume that  $P$ and $\omega$ are compatible. We have
  unearthed quite a bit of structure in the diagram \eqref{ref-7.1-36}.
  I have not seen this type of structure at other places.  
\begin{enumerate}
\item The diagram is $(3,2)$ periodic.
\item  Denote the
  horizontal and vertical maps by $h$ and $v$ respectively. 
$h$ and $v$ satisfy the following equations:
\begin{align*}
hv&=vh\\
v^2+\frac{1}{4}h^3&=1
\end{align*}
\item Every $2\times 2$ subdiagram is bicartesian.
\item The diagram is self dual up to a sign which must be applied to the
vertical maps. 
\end{enumerate}
\end{remark}
\begin{remark} 
\label{ref-7.6-52}
If we view the maps $\tilde{\omega}$ and $\tilde{P}$
as
endomorphisms of $D_A\oplus \Omega_A$ then in matrix form
they look like 
\[
\tilde{\omega}=\begin{pmatrix} 
\frac{1}{4} S &-\imath(P) \\
\imath(\omega) & 1
\end{pmatrix}
\qquad\qquad
\tilde{P}=\begin{pmatrix} 
1 &\imath(P) \\
-\imath(\omega) & \frac{1}{4} T
\end{pmatrix}
\]
We define a symmetric non-degenerate pairing between $D_A\oplus
\Omega_A$ and itself using the formula
\[
\langle (\delta,\eta),(\delta',\eta')\rangle=
\langle \delta,\eta'\rangle+\langle \delta',\eta\rangle
\]
For this pairing it is easy to check that $\tilde{\omega}$ is adjoint
to $\tilde{P}$. In other words $\tilde{\omega}$ is a unitary
transformation of $D_A\oplus \Omega_A$.  This suggest a connection
with a non-commutative version of Dirac geometry (yet to be created).
In the commutative case this connection is well understood \cite{BuCr}.
\end{remark}
Part (2) or Theorem \ref{ref-7.1-35} is proved in the same way as (1). It
remains to prove (3), i.e. the equivalence between the integrability 
conditions on $P$ and $\omega$. 

To do this we may follow a similar method as in \cite[Appendix]{VdB33}. We
start by expressing the integrability condition in terms of
Hamiltonian vector fields (see \S\ref{ref-2.3-3}). Throughout we assume that $P,\omega,\Phi$
satisfy ($\mathbb{P}$2)($\mathbb{P}$3)($\mathbb{B}$2)($\mathbb{B}$3)
as well as the compatibility condition \eqref{ref-7.3-32}. 
We let $\ldb-,-\rdb$ be the 
double bracket on $A$ induced by $P$ (see \eqref{ref-2.4-8}).
We recall that by Proposition \ref{ref-5.4-30}
 $P$ is quasi-Poisson
if and only if for all $a,b\in A$ we have the following identity in
$D_A\otimes A$
\begin{equation}
\label{ref-7.15-53}
\ldb H_a,H_b\rdb_l -H_{\ldb a,b\rdb'}(-)\otimes \ldb a,b\rdb''=
\frac{1}{4} [b,[a,E\otimes 1]_\ast]
\end{equation}
In order to work with Hamiltonian vector fields one needs some
formulas which are straightforward verifications
\begin{lemma} One has
\begin{align*}
i_{H_a}(\omega)&=da-\frac{1}{4}
[a,\Phi^{-1}d\Phi-d\Phi \Phi^{-1}]\\
i_{H_a}(d\Phi)&=-\frac{1}{2}[a,\Phi\otimes 1+1\otimes \Phi]_\ast
\end{align*}
\end{lemma}
\begin{proof} The first formula follows from the following computation
\begin{align*}
\imath(\omega)(H_a)&=\imath(\omega)\imath(P)(da)\\
&= (1-\frac{1}{4}T)(da)\\
&=da-\frac{1}{4}\sum_p [a,\Phi^{-1}d\Phi-d\Phi\Phi^{-1}]
\end{align*}
Now we prove the second formula
We
have
\[
i_{H_a}(d\Phi)=H_a(\Phi)=\ldb a,\Phi\rdb =-\ldb \Phi,a\rdb^\circ
\]
and according to \cite[Def 5.1.4]{VdB33}
\[
-\ldb \Phi,a\rdb^\circ=-\frac{1}{2} ((\Phi E+E\Phi)(a))^\circ                 
\]
We have
\[
(u E v)(a)=av \otimes u-v\otimes ua
\]
and thus
\[
(u E v)(a)^\circ=u\otimes av-ua\otimes v
\]
so that we get
\[
i_{H_a}(d\Phi)=-\frac{1}{2}(\Phi\otimes a- \Phi a\otimes 1+1\otimes a\Phi-a\otimes \Phi)
\]
as asserted.
\end{proof}
Using these formulas and a tedious computation one verifies the
following lemma
\begin{lemma} \label{ref-7.8-54} The formula \eqref{ref-7.15-53} always holds when evaluated 
on $d\Phi\otimes 1$. 
\end{lemma}
From this lemma we deduce
\begin{lemma} The element $P$ defines a quasi-Poisson algebra if and only if the
following identity in $\Omega_A\otimes A$ holds
\begin{equation}
\label{ref-7.16-55}
\imath_{\ldb H_a,H_b\rdb'_l}(\omega) \otimes \ldb H_a,H_b\rdb''_l
-\imath_{H_{\ldb a,b\rdb'}}(\omega)\otimes \ldb a,b\rdb''=
\frac{1}{4} [b,[a,\imath_E(\omega)\otimes 1]_\ast]
\end{equation}
\end{lemma}
\begin{proof}  By lemma \ref{ref-7.2-37} we have that
$
(\imath(\omega),e):D_A\r \Omega_A\oplus A(d\Phi)^\ast A
$
is injective. Hence to check that $\delta\in D_A$ is zero it is
sufficient to verify that $\imath_\delta(\omega)=0$ and
$e(\delta)=\delta(\Phi)'' (d\Phi)^\ast \delta(\Phi)'=0$. The latter holds
if $\delta(\Phi)=0$. The lemma now follows by writing \eqref{ref-7.15-53}
as $\delta'\otimes \delta''=0$ with $\delta'\in D_A$ and using Lemma 
\ref{ref-7.8-54}.
\end{proof}
To compute the part of \eqref{ref-7.16-55} involving $\ldb H_a,H_b\rdb_l $
we may use the formula \cite[(A.6)]{VdB33}.
\[
i_{H_a} \Lscr_{H_b}-\tau_{12} L_{H_b} \imath_{H_a}=
\imath_{\ldb H_a,H_b\rdb'_l}\otimes{\ldb H_a,H_b\rdb}''_l
+
\ldb H_a,H_b\rdb'_r \otimes \imath_{\ldb H_a,H_b\rdb''_r}
\]
Here $L_\delta=di_\delta+i_\delta d$ and $\Lscr_\delta={}^\circ L_\delta$. 
Applying the left hand side of this formula to $\omega$ we see that
everything is computable except the term $i_{H_a} \imath_{H_b} d\omega$
appearing as part of $i_{H_a}\Lscr_{H_b}(\omega)$. After a long and
tedious computation we arrive at the following lemma.
\begin{lemma} \label{ref-7.10-56} The formula \eqref{ref-7.16-55} holds if and only if
\[
\pr_1 i_{H_a} \imath_{H_b} d\omega=\frac{1}{6}\pr_1 i_{H_a} \imath_{H_b}(\Phi^{-1} d\Phi)^3
\]
holds for al $a,b$.
\end{lemma}
\begin{proof}[Proof of Theorem \ref{ref-7.1-35}(3)]
  It follows from lemma \ref{ref-7.10-56} that if $\omega$ is integrable then
  so is $P$.  To prove the converse
  we first do some computations
\begin{align*}
\imath_E d\omega&=-d(\imath_{E}\omega)\qquad \text{\cite[(A,7)]{VdB33}}\\
&=-\frac{1}{2}d(\Phi^{-1}d\Phi+d\Phi \cdot \Phi^{-1})\\
&=\frac{1}{2} (\Phi^{-1} d\Phi\cdot \Phi^{-1} d\Phi- d\Phi\cdot \Phi^{-1} d\Phi \cdot\Phi^{-1})
\end{align*}
On the other hand
\begin{align*}
\frac{1}{6}  \imath_{E}(\Phi^{-1} d\Phi)^3&=
\frac{1}{6}{}^\circ i_E(\Phi^{-1} d\Phi)^3\\
&=\frac{1}{6} {}^\circ (\Phi^{-1} (\Phi\otimes 1-1\otimes \Phi) \Phi^{-1} d\Phi\cdot
 \Phi^{-1} d\Phi+\cdots)\\
&=\frac{1}{2} (\Phi^{-1} d\Phi\cdot \Phi^{-1} d\Phi- d\Phi\cdot \Phi^{-1} d\Phi \cdot\Phi^{-1})
\end{align*}
So that we have deduced
\[
\imath_E d\omega= \frac{1}{6}  \imath_{E}(\Phi^{-1} d\Phi)^3
\]
Write
\[
\eta=d\omega-\frac{1}{6}(\Phi^{-1}d\Phi)^3
\]
Assuming that $P$ is integrable we need to prove that $\eta=0$ (modulo commutators).  We
already know by lemma \ref{ref-7.10-56} that $\pr_1 i_{H_a} \imath_{H_b}\eta=0$ and
by the computation above we have $\imath_E\eta=0$. 
 By \cite[Lemma A.5.2]{VdB33} we have $\imath(da)=H_a$.
  Hence by ($\mathbb{P}$3) we have $\sum_{a\in A}A H_a A+AEA=D_A$.
We observe
\begin{align*}
\pr_1 i_{E} \imath_{H_b} \eta&=-\pr_1 \tau_{12} i_{H_b} \imath_{E} \eta
\qquad \text{ \cite[(A.5)]{VdB33} } \\
&=0
\end{align*}
Hence by Prop.\
  \ref{ref-3.1.5-16}(2) it follows that $\imath_{H_b} \eta=0$.   Using $\imath_E\eta=0$ once again we conclude by Prop.\ \ref{ref-3.1.5-16}(1) that $\eta=0$ (modulo
commutators).
\end{proof}
\section{More general base rings and quivers}
\subsection{Generalities}
\label{ref-8.1-57}
In this section we show that the double quasi-Poisson brackets
constructed on (localized) path algebras of double quivers are
non-degenerate and hence Theorem \ref{ref-7.1-35} applies to them. 

However first we note that for quivers it is more
natural to use as base ring not a field but a direct sum of fields
indexed by the vertices. In \cite{CBEG} it was shown how to set up the
theory over an arbitrary semi-simple base ring $B$. In \cite{VdB33} we
worked over the base ring $B=ke_1+\cdots+ke_n$ with $e_p$ idempotent.
Since we rely on results from \cite{VdB33} we will do the same in this
section.

So assume that $B$ is as in the previous paragraph. We use
differentials and polyvector fields relative to $B$ (relevant
notations: $\Omega_{A/B}$, $\Omega_B A$, $D_{A/B}$, $D_B A$).  In this setting
the canonical element $E$ is defined as $\sum_p E_p $ where
\[
E_p(a)=ae_p\otimes e_p-e_p\otimes e_p a
\]
A (multiplicative) moment map $\Phi$ is now of the form $\sum_p
\Phi_p$ with $\Phi\in e_p A e_p $. With these conventions the
definitions and results in this paper go through verbatim. We will
accept this without further discussion.  For use below we introduce
the following convention. If $c=e_p A e_q$ then, changing standard
terminology, an \emph{inverse} of $c$ is an element $c^{-1}$ of $e_q A
e_p$ such that $c\cdot c^{-1}=e_p$, $c^{-1}c=e_q$. It is easy to see that
$c^{-1}$ is uniquely determined. 
\subsection{Fusion}
As in \cite{VdB33} our application to quivers depends on a process called
fusion (introduced in the commutative case in \cite{AKM}).  In the
case of quivers fusion amounts to gluing vertices but it is
beneficial to work somewhat generally.  We first construct and algebra $\bar{A}$ from $A$ by formally adjoining two
variables $e_{12}$, $e_{21}$ satisfying the usual matrix relations
$e_{uv}e_{wt}=\delta_{wv} e_{ut}$ (with $e_{ii}=e_i$).  The fusion
algebra of $A$ along $e_1,e_2$ is defined as
\[
A^f=\epsilon \bar{A}\epsilon
\]
where $\epsilon =1-e_2$. Clearly $\bar{A}$ is a $\bar{B}$-algebra and 
$A^f$ is a $B^f$-algebra. We will identify $B^f$ with 
$ke_1+ke_3+\cdots+ke_n$.

If $a\in A$ then we consider it as an element of $\bar{A}$  and
we write $a^f$ for $\epsilon a\epsilon+e_{12}a e_{21}$.  We extend
this convention to forms and polyvector fields.  Note that in \cite[\S
5.3]{VdB33} it was shown than the operations $\bar{(-)}$ and $(-)^f$ are
compatible with the formation of $D_BA$ and its Schouten bracket. A
similar result is true for $\Omega_B A$. The following result
was proved in \cite{VdB33}
\begin{theorems} \cite[Thm 5.3.1, 5.3.2]{VdB33} Assume that $(A,P,\Phi)$
  is a Hamiltonian double quasi-Poisson algebra (smooth as always).
  Then the same is true for $(A^f,P^{f\!\!f},\Phi^{f\!\!f})$ where
\[
P^{f\!\!f}=P^f-\frac{1}{2} E_1^f E_2^f
\qquad\text{and}\qquad
\Phi^{f\!\!f}_i=
\begin{cases}
\Phi^f_1\Phi^f_2&\text{if $i=1$}\\
\Phi^f_i&\text{if $i>2$}
\end{cases}
\]
\end{theorems}
In this section we prove the following result.
\begin{propositions}
\label{ref-8.2.2-58}
Assume that $(A,P,\Phi)$
  is a non-degenerate Hamiltonian double quasi-Poisson algebra. Then
the same is true for $(A^f,P^{f\!\!f},\Phi^{f\!\!f})$.
\end{propositions}
\begin{proof}
  In order to avoid confusing notations we define $F_p\in D_{B^f}
  (A^f)$ for $p\neq 2$ by $F_p(a)=ae_p\otimes e_p-e_p\otimes e_p a$.
It follows from \cite[(5.3)]{VdB33} that
\begin{equation}
\label{ref-8.1-59}
F_p=\begin{cases}
E_1^f+E_2^f&\text{if $p=1$}\\
E^f_p&\text{otherwise}
\end{cases}
\end{equation}
Since $A$ is non-degenerate the following map is surjective
\[
(\imath(P),c):\Omega_{A/B}\oplus \sum_p AE_pA\r D_{A/B}
\]
From this we deduce (by base extension) that the following map
is also surjective.
\[
(\imath(\bar{P}),c):\Omega_{\bar{A}/\bar{B}}\oplus \sum_p \bar{A}\bar{E}_p\bar{A}\r 
D_{\bar{A}/\bar{B}}
\]
By  \cite[Prop 4.2.1, lemma A.5.2]{VdB33} we have
\begin{equation}
\label{ref-8.2-60}
\imath(P)(db)=-\{ P,b\}
\end{equation}
Compatibility of fusion with Schouten brackets yields that
\[
(\imath(P^f),c):\Omega_{A^f/B^f}\oplus \sum_p \epsilon \bar{A}\bar{E}_p\bar{A}\epsilon \r
D_{A^f/B^f}
\]
is also surjective. In addition  we have $\epsilon \bar{A}\bar{E}_p\bar{A}\epsilon=
A^f E_p^f  A^f$.

It follows easily that 
\[
(\imath(P^{f\!\!f}),c):\Omega_{A^f/B^f}\oplus \sum_p A^f E^f_pA^f \r
D_{A^f/B^f}
\]
is surjective as well. Using \eqref{ref-8.1-59} we see that the current proposition
is proved provided we can show that $E^f_1$ is in the image of
$\imath(P^{f\!\!f})$, modulo $(F_p)_p$. We do this next. We have
\begin{align*}
\imath(P^{f\!\!f})(d\Phi^f_2)&=-\{P^f-\frac{1}{2}E_1^f E_2^f,\Phi_2^f\}
\end{align*}
We use some formulas we have already proved. I.e.\ the formula before
\cite[(5.5)]{VdB33} yields
\begin{align*}
\{P^f,\Phi_2^f\}&=-\frac{1}{2}(E_2^f \Phi_2^f +\Phi_2^f E_2^f )\\
&=-\frac{1}{2}((F_1-E_1^f) \Phi_2^f +\Phi_2^f (F_1-E_1^f ))\\
&=-\frac{1}{2}(F_1 \Phi_2^f +\Phi_2^f F_1)+\frac{1}{2}(E_1^f \Phi_2^f +
\Phi_2^f E_1^f )
\end{align*}
We also have (using the formulas after \cite[(5.5)]{VdB33})
\begin{align*}
\ldb E_1^f E_2^f,\Phi_2^f\rdb&=E_1^f\ast \ldb  E_2^f,\Phi_2^f\rdb-
\ldb E_1^f,\Phi_2^f\rdb\ast E_2^f\\
&=E_1^f\ast  (\Phi^f_2\otimes e_1-e_1\otimes \Phi^f_2)\\
&=\Phi^f_2\otimes E^f_1-E^f_1\otimes \Phi^f_2
\end{align*}
(``$\ast$'' represents the inner bimodule structure) so that we get
\[
\{ E_1^f E_2^f,\Phi_2^f\}=\Phi^f_2 E^f_1-E^f_1 \Phi^f_2
\]
Hence
\begin{align*}
  \{P^f-\frac{1}{2}E_1^f E_2^f,\Phi_2^f\}&=-\frac{1}{2}(F_1 \Phi_2^f
  +\Phi_2^f F_1)+\frac{1}{2}(E_1^f \Phi_2^f + \Phi_2^f E_1^f
  )-\frac{1}{2}(\Phi^f_2 E^f_1-E^f_1 \Phi^f_2)\\
&=-\frac{1}{2}(F_1 \Phi_2^f
  +\Phi_2^f F_1)+E_1^f \Phi_2^f
\end{align*}
So we are done.
\end{proof}
\subsection{Quivers}
Below we assume that $Q$ is a finite quiver whose vertices are indexed from
$1$ to $n$. We also use $Q$ to refer to the set of arrows of $Q$. The
head and tail of an arrow $a$ are denoted by $h(a)$, $t(a)$ respectively. 

Associated to $Q$ is the double quiver which has the same vertices as
$Q$ and arrows $\{a,a^\ast\mid a\in Q\}$ where $h(a^\ast)=t(a)$,
$t(a^\ast)=h(a)$. It will be convenient to write $(a^{\ast})^{\ast}=a$
and to define for $a\in \bar{Q}$, $\epsilon(a)=1$ if $a\in Q$ and
$\epsilon(a)=-1$ otherwise. We let $A$ be the path algebra of
$k\bar{Q}$ to which we adjoin the inverses of $1+aa^\ast$ for all
$a\in \bar{Q}$. 

For $a\in \bar{Q}$  one has a corresponding ``partial derivative'' in $D_{A/B}$
defined by
\[
\frac{\partial b}{\partial a}=
\begin{cases}
e_{t(a)}\otimes e_{h(a)} &\text{if $a=b$}\\
0&\text{otherwise}
\end{cases}
\]
It is easy to see that $D_{A/B}$ is a projective $A$-bimodule with generators
$\partial/\partial a\in e_{h(a)} D_{A/B}e_{t(a)}$. 
Note that $\partial/\partial a$ goes in the  ``opposite direction'' as $a$.

According to \cite[Thm 6.7]{VdB33} $A$ has a Hamiltonian double
quasi-Poisson structure given by
\[
P=\frac{1}{2}\left(\sum_{a\in \bar{Q}}
\left(\epsilon(a)
(1
+
a^\ast a)\frac{\partial\ }{\partial a} \frac{\partial\ }{\partial a^\ast}\right)
-
\sum_{a<b\in\bar{Q}}\left(\frac{\partial }{\partial  a^\ast} a^\ast-a\frac{\partial }{\partial  a}\right)\left(
\frac{\partial }{\partial  b^\ast} b^\ast-b\frac{\partial }{\partial  b}
\right)\right)
\]
\begin{equation}
\label{ref-8.3-61}
\Phi=\prod_{a\in \bar{Q}}(1+aa^\ast)^{\epsilon(a)}
\end{equation}
where ``$<$'' refers to an arbitrary ordering on the edges of $Q$, which is
also used to order the terms in the product \eqref{ref-8.3-61}.
\begin{propositions} This Hamiltonian double quasi-Poisson structure
  is non-degenerate.
\end{propositions}
\begin{proof}
This Hamiltonian double quasi-Poisson structure was constructed in 
\cite{VdB33} using fusion starting from (multiple copies of) the following
basic quiver. 
\begin{equation}
\label{ref-8.4-62}
\xymatrix{
1 \ar@/^/[rr]^a & &2\ar@/^/[ll]^{a^\ast}
}
\end{equation}
By Proposition\ \ref{ref-8.2.2-58} we may assume that $\bar{Q}$ is equal
to \eqref{ref-8.4-62}. The formula for $P$ then simplifies
to
\[
P=\frac{1}{2}\left((1+a^\ast a)\frac{\partial\ }{\partial a} \frac{\partial\ }{\partial a^\ast}
- (1+a a^\ast)\frac{\partial\ }{\partial a^\ast} \frac{\partial\ }{\partial a}\right)
\]
and by \cite[(6.2)]{VdB33}
\begin{align*}
E_1&=\frac{\partial\ }{\partial a^\ast}a^\ast -a\frac{\partial }{\partial a}\\
E_2&=\frac{\partial\ }{\partial a}a -a^\ast\frac{\partial }{\partial a^\ast}
\end{align*}
To check non-degeneracy we compute $\imath(P)(da)$ and $\imath(P)(da^\ast)$. 
We find
\begin{align*}
  \imath(P)(da)&=\frac{1}{2}{}^\circ\left((1+a^\ast a)i_{da}\frac{\partial\ }{\partial a} 
\frac{\partial\ }{\partial a^\ast}
+(1+a a^\ast)\frac{\partial\ }{\partial a^\ast} i_{da}\frac{\partial\ }{\partial a}
                    \right)\\
&=\frac{1}{2}\left(
\frac{\partial\ }{\partial a^\ast}(1+a^\ast a) 
+
(1+a a^\ast)\frac{\partial\ }{\partial a^\ast}
\right)
\end{align*}
A similar computation yields
\[
\imath(P)(da^\ast)=-\frac{1}{2}\left((1+a^\ast a)\frac{\partial\ }{\partial a}
+ \frac{\partial\ }{\partial a}   (1+a a^\ast)            \right)
\]
We have
\begin{align*}
a^\ast E_1+E_2 a^\ast&=-a^\ast a\frac{\partial\ }{\partial a}+ 
\frac{\partial\ }{\partial a}a a^\ast\\
&= -(1+a^\ast a)\frac{\partial\ }{\partial a}+\frac{\partial\ }{\partial a}(1+ a a^\ast )
\end{align*}
\begin{align*}
E_1a +a E_2 &=\frac{\partial\ }{\partial a^\ast}a^\ast a+ 
-a a^\ast\frac{\partial\ }{\partial a^\ast}\\
&=\frac{\partial\ }{\partial a^\ast} (1+a^\ast a)-(1+ a a^\ast )\frac{\partial\ }{\partial a^\ast}
\end{align*}
Hence
\begin{align*}
\imath(P)(da)&=(1+a a^\ast)\frac{\partial\ }{\partial a^\ast}\qquad \mod(E_1,E_2)\\
\imath(P)(da^\ast)&=-(1+a^\ast a)\frac{\partial\ }{\partial a} \qquad \mod(E_1,E_2)
\end{align*}
Since both $(1+aa^\ast)$ and $(1+a^\ast a)$ are invertible we are done.
\end{proof}
\def\cprime{$'$} \def\cprime{$'$}
\ifx\undefined\bysame
\newcommand{\bysame}{\leavevmode\hbox to3em{\hrulefill}\,}
\fi

\end{document}